\documentclass[11pt]{amsart}
\usepackage{amssymb,amsmath,amsthm}
\usepackage{mathrsfs,dsfont,a4wide}

\usepackage{latexsym}
\usepackage{amsfonts}
\usepackage{amsmath}
\usepackage{amssymb}
\usepackage{graphicx}
\numberwithin{equation}{section}
\newtheorem{theorem}{Theorem}[section]

\newtheorem{corollary}[theorem]{Corollary}
\newtheorem{definition}[theorem]{Definition}

\newtheorem{lemma}[theorem]{Lemma}

\newtheorem{remark}[theorem]{Remark}

\title
[Estimates for the dilatation of  $\sigma$-harmonic mappings]
{Estimates for  the dilatation of  $\sigma$-harmonic mappings}

\author[G. Alessandrini]
{Giovanni Alessandrini}
\address[G. Alessandrini]{Dipartimento di Matematica e Geoscienze, Universit\`a  di Trieste, Via Valerio 12/b, 34100 Trieste, Italia} 
\email[G. Alessandrini]{alessang@units.it}
\author[V. Nesi]
{Vincenzo Nesi}
\address[V.  Nesi]{Dipartimento di Matematica ``G. Castelnuovo", Sapienza, Universit\'a di Roma,
Piazzale A. Moro 2, 00185 Roma, Italy}
\email[V. Nesi]{nesi@mat.uniroma1.it}

\begin{document}

\begin{abstract}
We consider planar $\sigma$-harmonic mappings, that is mappings $U$ whose components $u^1$ and $u^2$ solve a divergence structure elliptic equation ${\rm div} (\sigma \nabla u^i)=0$, for $i=1,2$. We investigate whether a locally invertible $  \sigma$-harmonic mapping $U$ is also quasiconformal. Under mild regularity assumptions, only involving $\det \sigma$ and the antisymmetric part of $\sigma$, we prove quantitative bounds which imply quasiconformality.
\end{abstract}
\maketitle

\noindent{\small \textit{2000 AMS Mathematics Classification Numbers}: 30C62, 35J55}


\section{Introduction}
We study a certain class of homeomorphic solutions arising in several different context that we call $\sigma$-harmonic mappings.

\begin{definition}\label{def.H}
Consider an open, simply connected set $\Omega\subset
\mathbb R^2$. Given positive constants $\alpha$ and $\beta$, we
say that a measurable function $\sigma$, defined on $\Omega$ with
values into the space of $2\times2$ matrices,  belongs to the
class $\mathcal M (\alpha,\beta,\Omega)$  if one has
\begin{equation}\label{ellTM}
\begin{array}{ccrllll}
\sigma(x) \xi\cdot \xi&\geq& \alpha |\xi|^2&,&\hbox{for every $\xi \in \mathbb R^2$ and for a.e. $x\in\Omega$}\,,\\
\sigma^{-1}(x) \xi \cdot \xi &\geq & \beta^{-1}
|\xi|^2&,&\hbox{for every $\xi \in \mathbb R^2$ and for a.e.
$x\in\Omega$}\,.
\end{array}
\end{equation}
\end{definition}
\noindent
We say that $u$ is a $\sigma$-harmonic function if,
given $\sigma \in {\mathcal M}(\alpha, \beta,\Omega)$,   the function $u\in W^{1,2}_{\rm loc}(\Omega)$ is a weak solution to
\begin{equation}\label{u}
\begin{array}{ll}
{\rm div}(\sigma \nabla u )=0 & \hbox{in $\Omega$}\ .
\end{array}
\end{equation}
Here, and in what follows, we adopt the usual convention to identify $z=x_1 + i\,x_2\in\mathbb C$ with $x=(x_1,x_2)\in \mathbb R^2$.
A $\sigma$-harmonic map $U=(u^1,u^2):\Omega\to \mathbb R^2$ is simply a pair of $\sigma$-harmonic functions. If no conditions are imposed on the relationship between $u^1$ and $u^2$, the notion may degenerate to the uninteresting case when $u^1\equiv u^2$. We are interested to the case when $U$ is a homeomorphism and, from now on, we will always assume that it is sense preserving so that $\det DU\geq 0$ almost everywhere in $\Omega$ . The most classic example is the set of conformal mappings, that is holomorphic functions that are univalent, or, as is the same, injective. It is well know  that for holomorphic injective mappings $f$, the complex derivative $f^{\prime}$ does not vanish, or, in real notations, the Jacobian determinant does not vanish. Hence it is positive in $\Omega$ for conformal mappings, since they preserve the orientation. Another classical and  extensively studied example is the case of planar (univalent) harmonic mappings, \cite{duren}. 
 Lewy's celebrated theorem 
\cite{Lewy} extends to the harmonic planar mappings  the result which is valid for conformal mappings. So also for planar harmonic injective (sense preserving) mappings  $U$, one has $\det DU>0$  and $U$ is necessarily a diffeomorphism.

\noindent
For harmonic mappings, in particular for holomorphic one, $\sigma$ is the identity matrix, so no regularity issue arises. The map $U$ is real analytic. 
Another very classical case arises when $\sigma$ satisfies the additional constraint 
\begin{equation}\label{det=1}
\det \sigma=1\ .
\end{equation}
Indeed, let us recall that 
given $\sigma \in {\mathcal M}(\alpha, \beta,\Omega)$   and  a weak solution $u\in W^{1,2}_{\rm loc}(\Omega)$  to
\eqref{u},
then there exists $\tilde u \in W^{1,2}_{\rm loc}(\Omega)$, called the \emph{stream function} of $u$, such that, setting
\begin{equation}\label{J}
 J:=\left(
\begin{array}{ccc}
0&-1\\
1&0
\end{array}
\right)\,,
\end{equation}
\noindent
 one has
\begin{equation}\label{e2.3}
\begin{array}{lllll}
 \nabla \tilde u= J \sigma \nabla u &\hbox{in}&\Omega\,.
 \end{array}
\end{equation}
If, in addition, $\sigma$ is symmetric and \eqref{det=1} holds, then one easily verifies that $\tilde{u}$ also satisfies \eqref{u}. Hence $U=(u,\tilde{u})$ is $\sigma$-harmonic.

\noindent
For any solution to \eqref{u},
setting
\begin{equation}
\label{e2.4}
F=u+i \tilde u\,,
\end{equation}
one has
$F=u+i \tilde u \in W^{1,2}_{\rm loc}(\Omega;\mathbb R^2)$ and one writes,  in complex notations,
\begin{equation}\label{1storder}
\begin{array}{ll}
F_{\bar{z}}=\mu F_z +\nu \bar{F_z}\ & \hbox{in $\Omega$}\ ,
\end{array}
\end{equation}
where, the so called complex dilatations $\mu , \nu$ are given by
\begin{equation}\label{SNU}
\begin{array}{llll}
\mu=\frac{\sigma_{22}-\sigma_{11}-i(\sigma_{12}+\sigma_{21})}{1+{\rm
Tr\,}\sigma +\det \sigma}& \ ,&\nu =\frac{1-\det \sigma
+i(\sigma_{12}-\sigma_{21})}{1+{\rm Tr\,}\sigma +\det \sigma}\ ,
\end{array}
\end{equation}
and satisfy the following
ellipticity condition
\begin{equation}\label{ellQC}
|\mu|+|\nu|\leq \frac{K-1}{K+1} \,,
\end{equation}
 for some $K\geq 1$ only depending on
$\alpha, \beta$. 
As is well known, \eqref{e2.4}-\eqref{ellQC}, imply that 
\begin{equation}\label{dist1}
|DF|^2\leq \left(K+\frac 1 K\right) \det DF\,,
\end{equation}
where $|B|=({\rm Trace}(B B^T))^{1/2}$ denotes the Hilbert-Schmidt norm of the matrix $B$. Any $W^{1,2}_{\rm loc}(\Omega)$ mapping satisfying \eqref{dist1} almost everywhere, is called $K$-quasiregular.
In fact a classical equivalent way of introducing  quasi regular mappings is in terms of ACL (Absolutely Continuous on almost very Line) orientation preserving mappings for which there exist $K\geq 1$ such that \eqref{dist1} holds.

\noindent
Associated to the map $F$ we may associate various types of so-called distortions or dilatations, which turn out to be equivalent.
Let us consider the following one
\begin{equation}\label{dist2}
\begin{array}{ccc}
d(F)=\frac{|DF|^2}{2 \det DF} \,.
\end{array}
\end{equation}
\noindent
The principal object of the present paper is to establish conditions on $\sigma$ under which a $\sigma-$harmonic mapping $U$ is quasiregular, that is $d(U)$ is bounded.
\noindent
For reasons that shall be made transparent later on, it will be convenient to introduce a slight variation of the dilatation function and introduce a new but equivalent dilatation $d^{\sigma}$ as follows. For a sense preserving mapping $U$, that is, for a mapping such that $\det DU\geq0$, we define
\begin{equation}\label{dsigma}
d^{\sigma}(U):=\frac{{\rm Trace} (DU \sigma DU^T)}{2 \det DU}=\frac{\sigma \nabla u^1 \cdot \nabla u^1+\sigma \nabla u^2 \cdot \nabla u^2}{2 \det DU}\,.
\end{equation}
In view of the ellipticity of the matrix $\sigma$, it is clear that $d$ and $d^{\sigma}$ are uniformly bounded one from each other by positive constants.
It must be emphasized that, in general, the second component of a $\sigma$-harmonic mapping need not be the stream function of the first one, already at the level of harmonic mappings. So $\sigma$-harmonic mappings are indeed a much wider class of the set of quasiconformal mappings.

\noindent A very natural question is whether the dilatation $d^{\sigma}(U)$ may blow up for general $\sigma$-harmonic maps.

\noindent In fact,  already in the case of harmonic mappings, i.e. when $\sigma$ is the identity matrix $I$, one may have $d^I(U)$ diverging. To see this, it suffices to consider, in complex notation, the homeomorphism 
\begin{equation}
U(z)=z+\frac 1 2 \bar{z}^2
\end{equation} 
mapping the closed unit disk $B(0,1)$ univalently onto a closed bounded set the boundary of which is an hypocycloid. One easily checks the Jacobian determinant vanishes identically on $\partial B(0,1)$ and hence the dilatation diverges at any point on it. However, by smoothness and Lewy's Theorem, the dilatation may only blow up at the boundary.  

\noindent
For $\sigma$-harmonic mappings, instead, the situation is dramatically different. In fact, not only the dilatation may blow up at interior points but, actually, this is the generic situation.

\noindent
A typical example of $\sigma$ for which one can construct a univalent $\sigma$ harmonic mapping the dilatation of which diverges, is given in \cite{ANJAM}, Example 3.1. A much wider class is given in the same paper, Example 3.2. In the same paper it is proved, roughly speaking, that if for a given $\sigma$ there exists a univalent mapping $U=(u^1,u^2)$ which is a solution to
\begin{equation}\label{u1u2}
\begin{array}{ll}
{\rm div}(\sigma \nabla u^i )=0 & \hbox{in $\Omega$ for $i=1,2$}\,,
\end{array}
\end{equation}
that is such that $d^{\sigma}(f)$ is locally bounded away from the boundary of $\Omega$, then {\it all the univalent mappings $U$ which are solutions to \eqref{u1u2}} are quasiconformal, see also \cite{nf}.
The class of all such nice $\sigma$'s, denoted by $\Sigma_{\rm qc}$ is rather ``small", that is nongeneric. More precisely, $\Sigma_{\rm qc}$ is a set of first Baire category as proved in Theorem~4.1 of \cite{ANJAM}. The identity matrix is one of such nice $\sigma$'s. So, for given $\sigma$, in order to have $\sigma\in \Sigma_{\rm qc}$, it suffices to find a single $\sigma$-harmonic mappings with bounded dilatation.

\noindent
One may wonder how to exhibit large classes of univalent $\sigma$-harmonic mappings. The first result for the case of $\sigma$ not equal to the identity, was established in \cite{bmn}, Theorem~3.1, following previous work by Bauman and Phillips \cite{BP}. The idea, going back to Kneser \cite{k} and Rad\`o \cite{r}, is to deduce the univalence by some informations about the Dirichlet boundary data. In \cite{bmn} one assumes that the boundary data is the restriction to the boundary of a smooth sense preserving diffeomorphism $\Psi=(\Psi^1,\Psi^2)$ of an  open, simply connected set $\Omega$ with smooth boundary, onto a convex set. Under these assumptions one proves that $U$ is itself a self preserving homeomorphism so that $\det DU\geq 0$ almost everywhere. Under more stringent regularity assumptions on $\partial \Omega$, $\Psi$ and $\sigma$, one deduces that $DU$ is smooth and $\det DU>0$.

\noindent
A deeper result, in this direction, was proved in \cite{ANARMA}, Theorem~4. In that paper a Kneser-Rad\`o type of theorem is proved under essentially minimal hypotheses. For the kind of questions we address in the present paper, however, the main result is  a different one. To state it in a a precise form, we recall that, given an  open set $D\subset \mathbb R^2$, $\phi \in L^1_{\rm loc} (D)$ belongs to  ${\rm BMO} (D)$ if
\begin{equation*}
\| \phi \| _{*} = \sup_{Q\subset D} \left( \frac{1}{\mid Q\mid}
\int_Q \mid \phi - \phi_Q\mid \right) <\infty\,,
\end{equation*}
where $Q$ is any square in $D$ and $\phi_Q= \frac{1}{\mid Q\mid} \int_Q \phi$. 

\noindent
Recall also that the normed space
 $({\rm BMO} (D), \|\cdot\|_*)$ is in fact a Banach space.
  \begin{theorem}
\label{t3.1}
Let $\Omega$ be an open subset of $\mathbb R^2$, let $\sigma \in {\mathcal M}(\alpha,\beta,\Omega)$  and
let $U \in W^{1,2}_{\rm loc}(\Omega,\mathbb R^2)$ be a $\sigma$-harmonic mapping  which is locally one-to-one and sense preserving. For every $D\subset \subset \Omega$, we have
\begin{equation}
\label{e3.1}
\begin{array}{ll}
\log (\det DU) \in {\rm BMO} (D)\,.
\end{array}
\end{equation}
\end{theorem}
 \noindent
 One obtains by \eqref{e3.1}, that the Jacobian determinant has a negative exponent of integrability and, in particular, that  $\det DU>0$ almost everywhere. This theorem, which relies on deep results due to Reimann \cite{re}, Bauman \cite{bau} and Fabes and Strook \cite{fs}, was proved in \cite{ANARMA}, Theorem~5, in the case when $\sigma$ is symmetric. The extension to the nonsymmetric case elaborates on results in \cite{nf} and can be found in \cite{ANFINNICO}, Theorem~3.1.  In fact, in \cite{ANFINNICO} Theorem~4.1 and Remark~4.3,
an even stronger result is obtained, that is
\begin{equation}\label{peso}
\det DU \in A_{\infty}(D),
\end{equation}
where $A_{\infty}$ is the class of Muckenhoupt weights \cite{cf}. Indeed, in Theorem~\ref{theorem2} below, we shall make use of this result.

\noindent
Applications of these achievements may be found in several context, ranging form establishing sharp bounds for effective conductivity as in \cite{achn} and \cite{acn}, or in establishing a conjecture in the field of quasiconfomal mappings about the stability of Beltrami systems under $H$-convergence,  as in \cite{ANFINNICO}. Other applications have been given in inverse problems, see, for instance,  the review article by \cite{balinside}.

\noindent
In the present paper 
we investigate sufficient conditions on $\sigma$ such that $d^{\sigma}(U)$ is locally bounded. In \cite{ANJAM}, Theorem~2.2, we proved that if $\sigma$ is symmetric and $\det \sigma$ is H\"older continuous, then $d^{\sigma}(U)$ is locally bounded, however no quantitative estimate was obtained. Here, assuming that the antisymmetric component $\sigma_{12}-\sigma_{21}$ and the determinant $\det \sigma$ are Lipschitz continuous, we achieve a concrete local upper bound on $d^{\sigma}$ in the form of a Harnack type inequality in Theorem~\ref{theorem1}.

\noindent
Next, since we know from \cite{ANJAM}, as already remarked, that in order to have $\sigma\in \Sigma_{qc}$, it suffices to exhibit a single $\sigma$-harmonic mapping whose dilatation is bounded, we focus on the special class of so-called periodic $\sigma$-harmonic mappings, which are of special relevance in homogenization and which have beee analyzed in depth in \cite{ANARMA} and in \cite{ANCOCV}.
More precisely, denoting by $W^{1,2}_{\sharp}(\mathbb R^2,\mathbb R^2)$ the space of $W^{1,2}_{\rm loc}$ mappings whose components are $1$-periodic in each variable, and assuming that $\sigma \in \mathcal M (\alpha,\beta,\Omega)$ is also $1$-periodic, in each variable, we consider $U=(u^1,u^2)\in W^{1,2}_{\rm loc}(\mathbb R^2,\mathbb R^2)$ such that
\begin{equation}\label{Uper}
\left\{
\begin{array}{lllll}
{\rm div}(\sigma \nabla u^i )= 0&{\rm in}&\mathbb R^2\,,\\ 
U-x\in W^{1,2}_{\sharp}(\mathbb R^2,\mathbb R^2)&&\,.
\end{array}
\right.
\end{equation}
For the solution to \eqref{Uper},  we prove a global uniform bound on $d^{\sigma}(U)$ in the final Theorem~\ref{theorem2}.

\section{An elliptic equation}
\noindent
Throughout the present Section, we assume that $U$ is a locally univalent,  sense preserving $\sigma$-harmonic mapping $U=(u^1,u^2)$ in $\Omega$. The main motivation for the introduction of the quantity $d^{\sigma}$ comes for the following decomposition.
Denote
\begin{equation}\label{wi}
\begin{array}{ll}
w^i=\frac{\det DU}{\sigma \nabla u^i\cdot \nabla u^i}\,,& \hbox{for $i=1,2$}\,.
\end{array}
\end{equation}
Obviously we have
\begin{equation}\label{wvsd}
d^{\sigma}(U)=\frac 1 2 \left(\frac{1}{w^1}+ \frac{1}{w^2}\right).
\end{equation}
We shall show in the next Lemma that, under few hypotheses on $\sigma$, each $w^i$ solves a divergence structure elliptic equation. In fact is is well known that if $\sigma$ is the identity, then indeed $w^1$ and $w^2$ are harmonic \cite{duren}.

\begin{lemma}\label{lemma1}
Assume, for some $E\in (0,+\infty)$,

\begin{equation}\label{Lip}
||\nabla\det \sigma||_{L^{\infty}(\Omega)}+
||\nabla (\sigma_{12}-\sigma_{21})||_{L^{\infty}(\Omega)}\leq E.
\end{equation}
Then, for any $i=1,2$, $w^i\in W^{1,2}_{\rm loc}(\Omega)$ and it is a weak solution to

\begin{equation}\label{neweq}
\begin{array}{llll}
{\rm div}(\sigma \nabla w^i +w^i B^i)=0&\hbox{\rm{in} $\Omega$}\ ,
\end{array}
\end{equation}
where
\begin{equation}
\begin{array}{lll}
B^i=
\left\{
\begin{array}{lcc}
 \displaystyle{\frac{1}{\sigma \nabla u^i \cdot\nabla u^i}}
\left[ \left(J\nabla u^i \cdot\nabla c\right) 
  J^T \nabla u^i + \left(J\nabla u^i \cdot\nabla b\right) 
\sigma\nabla u^i 
  \right]\,,
  &\hbox{\rm{where}}&\nabla u^i\neq 0
 \\
\\
0\,,&\hbox{\rm{where}}&\nabla u^i= 0
\end{array}
\right.
\end{array}
\end{equation}
and
\begin{equation}
\begin{array}{lll}
b=\sigma_{12}-\sigma_{21}&,& c=\det \sigma\,.
\end{array}
\end{equation}
\end{lemma}
\begin{remark}
Being $U$ locally univalent, we know that $\det DU$ may vanish only on sets of zero measure. Even more so, $|\nabla u^i|$, may vanish only on sets of zero measure. See \cite{AJMAA} for a general discussion of unique continuation for two-dimensional elliptic equations.
Observe that
\begin{equation*}
||B_i||_{L^{\infty}(\Omega)}\leq C_0\,E,
\end{equation*}
where $C_0$ only depends on $\alpha$ and $\beta$. Finally note that $w^i$ are actually $\sigma$-harmonic when $\det \sigma$ and the antisymmetric part of $\sigma$ are constant. 
\end{remark}
\noindent
{\bf Proof of Lemma~\ref{lemma1}.} It suffices to provide a formal proof of \eqref{neweq}. Standard smoothing arguments and interior a-priori estimates will guarantee that $w^i\in W^{1,2}_{\rm loc}$ and the validity of \eqref{neweq} with no further assumptions other than \eqref{Lip}.

\noindent
We concentrate on $w^1$, the same arguments shall apply to $w^2$.
It was shown in \cite{ANARMA} and \cite{ANFINNICO} that,  being $U$ injective, setting $f=u^1+ i\,\tilde{u}^1$, $f$ is locally quasiconformal and, introducing
\begin{equation*}
V=U\circ f^{-1}\,,
\end{equation*}
one has
\begin{equation*}
w:=\det DV=\frac{\det DU}{\det DF} \circ f^{-1}=w^1\circ f^{-1}\,.
\end{equation*}
Furthermore, see \cite{ANFINNICO}, $w$ is a distributional solution to the following elliptic equation in adjoint form
\begin{equation}\label{adeq}
w_{x_1x_1}+(b\circ f^{-1}\,w)_{x_1x_2}+(c\circ f^{-1}\,w)_{x_2x_2}=0\,.
\end{equation}
From now on we will simply write, with a slight abuse of notation $b$ and $c$ in \eqref{adeq} instead of
\begin{equation*}
b\circ f^{-1}\hskip0.5cm\hbox{and}\hskip0.5cm c\circ f^{-1}\,.
\end{equation*}
In view of the bound \eqref{Lip}, $b$ and $c$ are weakly differentiable, hence  \eqref{adeq} can be rewritten as a divergence structure equation
\begin{equation}\label{adeqdiv}
w_{x_1x_1}+(b \,w_{x_2} +b_{x_2} \,w)_{x_1}+(c\,w_{x_2} +c_{x_2}\,w)_{x_2}=0\,.
\end{equation}

\noindent
Pulling back through the quasiconformal mapping $f,$ we end up with equation
\eqref{neweq}.
\section{Bounds on the dilatation}
In what follows we are interested in pointwise bounds on the dilatation function $d^{\sigma}(U)$ for a given $\sigma$-harmonic mapping $U$. Therefore, in order to emphasize the dependence on the variable $x\in \Omega$, we shall write for short
\begin{equation}
d^{\sigma}(x)\hskip0.5cm\hbox{instead of}\hskip0.5cm d^{\sigma}(U)(x)\,.
\end{equation}
Also, for the sake of simplicity, but with no loss of generality, we shall assume that the ellipticity constant in \eqref{ellTM} satisfy
\begin{equation}\label{K}
\alpha^{-1}=\beta=K\,.
\end{equation}

\begin{theorem}\label{theorem1}
Let $\sigma$ satisfy \eqref{ellTM} with the convention \eqref{K} and assume that  $U:\Omega\to\mathbb R^2$ be a locally univalent and sense preserving $\sigma$-harmonic mapping. Assume, in addition \eqref{Lip}.
Then $d^{\sigma}<+\infty$ in $\Omega$, it is H\"older continuous and for any $A\subset \subset \Omega$, there exists $H>0$ only depending on $K, E, A$ and $\Omega$, such that the following Harnack's inequality holds

\begin{equation}\label{Harnack}
d^{\sigma}(x)\leq H \,d^{\sigma}(x^{\prime}) \,\,\hbox{for all $x, x^{\prime}\in A$}\,.
\end{equation}
\end{theorem}
\noindent{\bf  Proof}. Recall \eqref{wvsd}. In view of Lemma~\ref{lemma1}, by De Giorgi-Nash-Moser theory \cite{gt}, the thesis follows.

\vskip.1cm
\noindent
As a straightforward consequence, we obtain the following
\begin{corollary}\label{corollary} Under the same assumptions as in Theorem~\ref{theorem1}, we have
\begin{equation}\label{stima1}
||d^{\sigma}||_{L^{\infty}(A)}\leq  H \left(\frac{1}{|A|}\int_A(d^{\sigma})^{\delta}\right)^{\frac{1}{\delta}}
\end{equation}
for any positive $\delta$.
\end{corollary}

\begin{remark}
As already pointed out, under the assumptions of Theorem~\ref{theorem1}, the exists a positive $\epsilon$ such that $(\det DU)^{-\epsilon}$ in integrable. Therefore, by Schwartz inequality, one can check  that the right hand side of \eqref{stima1} is bounded for some small $\delta>0$ even when the condition \eqref{Lip} is not assumed.
\end{remark}

\begin{theorem}\label{theorem2}
Assume, in addition to the hypotheses of Theorem~\ref{theorem1}, that $\sigma$ is 1-periodic in each variable on the whole $\mathbb R^2$. Consider the solution to the periodic problem \eqref{Uper} with ``unit cell" $Q_0=[0,1]\times[0,1]$.
Then there exists $M>0$ only depending on $K$ and $E$ such that
\begin{equation}
||d^{\sigma}||_{L^{\infty}(Q_0)}\leq M\,.
\end{equation}
\end{theorem}
\noindent{\bf Proof.}
The solution $U:\mathbb R^2\to\mathbb R^2$ is globally univalent, as proved in \cite{ANARMA}, Theorem~1. By Theorem~4.1 in \cite{ANFINNICO} and by the  theory of Muckenhoupt weights, \cite{cf}, we can find $p>1$ and $C>0$ only depending on $K$ such that

\begin{equation}\label{disC}
\frac{1}{|R|} \int_R \det DU dx \left(\frac{1}{|R|} \int_R (\det DU)^{-\frac{1}{p-1}}dx \right)^{p-1}\leq C\,,
\end{equation}
for every square $R\subset \mathbb R^2$. By Corollary~\ref{corollary}, one has
\begin{equation}
||d^{\sigma}||_{L^{\infty}(Q_0)}\leq  \frac{H}{2} \left(\int_{Q_0} (|DU|_{\sigma})^{2\delta} (\det DU)^{-\delta}dx\right)^{\frac{1}{\delta}},
\end{equation}
where, for given $\sigma\in \mathcal M (\alpha,\beta,\Omega)$, we set $|B|_{\sigma}:=({\rm Trace}(B \sigma B^t))^{1/2}$, a Hilbert-Schmidt norm of the matrix $B$ ``weighted'' with that matrix $\sigma$.

\noindent
By Schwartz's inequality
\begin{equation*}
||d^{\sigma}||_{L^{\infty}(Q_0)}\leq  \frac{H}{2} \left(\int_{Q_0}(|DU|_{\sigma})^{4\delta} dx\right)^{\frac{1}{2\delta}}
 \left(\int_{Q_0}(\det DU)^{-2 \delta}dx\right)^{\frac{1}{2\delta}}.
  \end{equation*}
 Picking $0<\delta<\min\left(\frac 1 2, \frac{1}{2(p-1)}\right)$ and recalling that $|Q_0|=1$, we deduce by H\"older's inequality that
\begin{equation*}
\int_{Q_0}(|DU|_{\sigma})^{4\delta}dx\leq
 \left(\int_{Q_0}(|DU|_{\sigma})^{4\delta s}dx\right)^{\frac{1}{s}},\,\,\forall s>1\,.
  \end{equation*}
  Hence, picking $s=\frac{1}{2\delta}$, we have
  \begin{equation*}
\left(\int_{Q_0}(|DU|_{\sigma})^{4\delta}dx\right)^{\frac{1}{2\delta}}\leq
\int_{Q_0}(|DU|_{\sigma})^{2}dx\,.
  \end{equation*}
  On the other hand, using again H\"older's inequality and the fact that $|Q_0|=1$, we have
  \begin{equation*}
\int_{Q_0}(\det DU)^{-2\delta}dx\leq
 \left(\int_{Q_0}(\det DU)^{-2\delta r}dx\right)^{\frac{1}{r}},\,\,\forall r>1
  \end{equation*}
  and choosing $r=\frac{1}{2\delta (p-1)}$, we have, by \eqref{disC}
  \begin{equation*}
\left(\int_{Q_0}(\det DU)^{-2\delta}dx\right)^{\frac{1}{2\delta}}\leq
 \left(\int_{Q_0}(\det DU)^{-\frac{1}{p-1}}dx\right)^{p-1}\leq
 C \left(\int_{Q_0}(\det DU)dx\right)^{-1}.
  \end{equation*}
 By standard energy estimates
   \begin{equation*}
\int_{Q_0}(|DU|_{\sigma})^{2}dx\leq  \int_{Q_0}{\rm Trace}\,(\sigma) dx\leq 2K
  \end{equation*}
  and, since $\det DU$ is a null-lagrangian, by the standard Area theorem, one has
   \begin{equation*}
\int_{Q_0}\det DU dx=\int_{Q_0}\det I dx=1.
  \end{equation*}
  Hence
  \begin{equation*}
||d^{\sigma}||_{L^{\infty}(Q_0)}\leq  C H K\,.
\end{equation*}
\vskip0.5cm
\noindent
{\bf Acknowledgements} G.A. was supported by 
Universit\`a degli Studi di Trieste FRA 2012, Problemi Inversi.
V.N. was supported by PRIN Project 2010-2011 ``Calcolo delle Variazioni''.

\end{document}